\documentclass[12pt]{article}
\usepackage{amssymb}
\usepackage{theorem}

\newtheorem{theorem}{Theorem}

\newtheorem{lemma}[theorem]{Lemma}
{\theorembodyfont{\rmfamily}

}
\include{epsf}
\begin{document}

\def\CP{\mathbb{C}{\rm P}}
\def\tr{{\rm tr}\,}
\def\endproof{{$\Box$}}

\title{Foliations with unbounded deviation on $T^2$}

\author{Dmitri Panov}

\maketitle

\section{Introduction}
The goal of this note is to prove the following statement.

\begin{theorem} There exists a $C^{\infty}$ foliation with 3 singular 
points on the two-dimensional torus such that any lifting of a leaf of this 
foliation on  the universal covering of the torus is a dense subset of 
the covering. 
\end{theorem}

This foliation is nonorientable and it has a transversal measure. 
The measure can be given locally out of the singularities by a closed
1-form.

{\bf Half translation surfaces.} Let $M$ be a two-dimensional surface with a
flat metric  which has only conical singularities. The surface is called 
a {\it half translation} if the holonomy along any path on the surface is 
either $Id$ or -$Id$. It is clear that the angles around the conical 
points should be an integer multiple of $\pi$. The surface is called 
a {\it translation} if the holonomy is always Id.

A half translation structure gives rise to a family of 
foliations on the surface parametrized by the circle $S^1$. 
To construct a foliation from this family fix a tangent vector on the 
surface and move it to all the points of the surface by means of 
parallel translation. Since the holonomy is $\pm Id$ we obtain a foliation.

A half translation surface can have only a finite group of
automorphisms (with the only one exception of a flat torus).
But if we forget about the metric and remember only the underlying
affine structure on the surface we can have an interesting group of 
automorphisms preserving affine structure.  

The idea of the proof of the theorem 1 is the following.
We construct a half translation torus with 3 singular points such 
that the underlying affine structure has a big group of
automorphisms. In particular there is an automorphism of
the torus  which acts trivially  in the first homology of the torus
but is a locally hyperbolic map (in affine coordinates on the
torus). The expanding foliation of this automorphism will be
the foliation mentioned in the theorem 1.

{\bf Remarks}. There is a classical resalt of A.Weil [We] about  
oriented foliation with  finite number of singularities on $T^2$.
The {\it deviation} of such a foliation is always bounded i.e.  
any lifting of any leaf of the foliation to the universal covering of 
the torus is contained in a finite neighborhood of a strait line. 
First example of a nonoriented foliation
with unbounded deviation was constructed by A.D. Anosov [An1], [An2].
In [EMZ] you can find a good introduction to translation surfaces.

\section{Construction of the foliation}

{\bf Construction of the half translation structure on the torus.}
The  half translation structure on the torus  will be obtained from a usual 
flat torus by the operation of folding. 

{\bf Definition of  folding.} Let $[AB]$ be a geodesic segment on a 
flat surface M. Cut M along $[AB]$ and denote by $C_+$ and $C_-$ 
two centers of different shores of the cut. 
Glue together  segment $[AC_+]$ with $[C_+B]$ and $[AC_-]$ with $[C_-B]$.
On the resulting surface we obtain two new singular points with 'angles' $\pi$
(these points appear from $C_+$ and $C_-$) and a singular point with 
'angle' $4\pi$. 

\begin{center}
\
\epsfbox{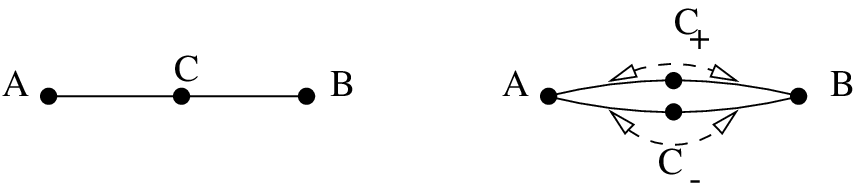}
\end{center}

The torus we are looking for is the flat torus
$\mathbb{R}^2 /(3\mathbb{Z}\oplus \mathbb{Z})$ folded along the segment 
$[(-1,0);(1,0)]$. We denote this torus by $T$. The structure of the folded
torus is given on the fig.2. The arrows connect pieces of the border of 
the rectangle that should be glued. 
\begin{center}
\
\epsfbox{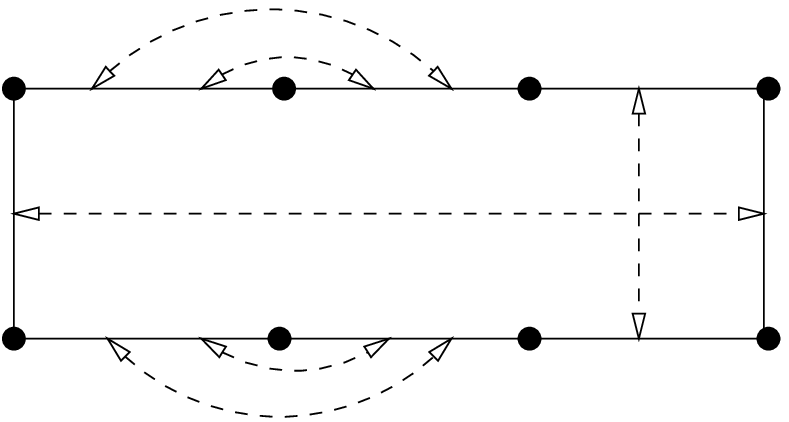}
\end{center}

{\bf Affine automorphisms of $T$. Decomposition on cylinders.}   
Consider the $\it horizontal$ foliation (i.e. the foliation
whose leaves are parallel to the $x$-axis) on the torus $T$  and 
remove the leaf passing through singular points. It cuts $T$
into the cylinder with height 1 and length 3. The horizontal foliation 
fills the cylinder by closed leaves of length 3. There is an obvious
automorphism $\gamma_h$ of  $T$ preserving the horizontal foliation. It
fixes the singular leaf and it makes a twist of the cylinder.
This map can be given by the formula  $\gamma_{h}(x,y)=(x+3y(mod\ 3),y)$. 

Now consider the $\it vertical$ foliation. This foliation has three 
connected singular leaves. They split $T$ into two cylinders
$1\times 1$ and $2\times 1$ (here and later the cylinder $m\times n$ is
the cylinder with the height $m$ and the length $n$). There is an 
automorphism $\gamma_{v}$ of $T$ which  makes the full twist of the 
cylinder $1\times 1$ and it makes the half  twist of the 
cylinder $2\times 1$. The automorphism 
$\gamma_{v}$ fixes the singular point having the conical angle $4\pi$ 
and permutes two singular points having the conical angle $\pi$.

There is a natural map from the group of affine automorphisms of $T$ to
$PSL(2,\mathbb{Z})$. It is given by the action of the automorphisms on
the tangent space of $T$.
The images of $\gamma_h$ and $\gamma_v$ are given by 
$$\left(\begin{array}{cc}
1&3\\
0&1\\
\end {array}\right) ,                              
\left(\begin{array}{cc}
1&0\\
1&1\\
\end {array}\right)$$                      
These matrixes are defined only up to a sign,   
                               
The automorphism we are looking for is 
$(\gamma_{v}\gamma_{h}^{-1}\gamma_v)^4$.

{\bf Construction of the foliation.}
We denote by $\mathcal F_1$ the expanding foliation of the map
$\gamma_{v}\gamma_{h}^{-1}\gamma_v$
and by $\mathcal F_2$ the contracting foliation. 
These foliations are parallel to 
the directions $(\sqrt{3},1)$ and $(\sqrt{3},-1)$. 

\begin{lemma} Any leaf of the foliation $\mathcal F_1$ is a dense
 subsets  of the torus.
\end{lemma}

The proof of this lemma is standard. It follows from the absence of 
separatrices joining singular points of the foliation.

\section{Proof of the theorem}
{\bf Action on homology.}  
We note by $e_1$ the class in the homology group of a nonsingular 
leaf of the horizontal foliation and we note by $e_2$ the class of 
homology of a leaf of the vertical cylinder $1\times 1$. It is clear that 
${e_1,e_2}$ form a basis for $H_1(T,\mathbb{Z})$.

\begin{lemma} The action of the induced homomorphism  
$\gamma_{v}\gamma_{h}^{-1}\gamma_v$  on the fist homology of $T$ 
with respect to the basis $(e_1,e_2)$ is given by the matrix

$$\left(\begin{array}{cc}
0&-1\\
1&0\\
\end {array}\right)$$

\end{lemma}
Proof. Since $\gamma_h$ makes the single twist of the horizontal cylinder
$1\times 3$ and $\gamma_v$ makes the single twist of the vertical 
cylinder $1\times 1$ the action of $\gamma_h$ and $\gamma_v$ in 
$H_1(T,\mathbb{Z})$ is given by the matrixes 

$$\left(\begin{array}{cc}
1&1\\
0&1\\
\end {array}\right)
\left(\begin{array}{cc}
1&0\\
1&1\\
\end {array}\right)$$

It proves the lemma.\hfill $\square$

Lemma 3 implies that the induced action of 
$(\gamma_{v}\gamma_{h}^{-1}\gamma_v)^4$ on $H_1(T,\mathbb{Z})$ is trivial.

{\bf Action on the universal covering of torus $T$.} 
The universal covering of $T$ can be viewed as the usual plane
$\mathbb{R}^2$ folded
along the segments $[(n,3k-1), (n,3k+1)]$, $n,k\in \mathbb{Z}$. 
After folding the origin $(0,0)$ of the plane gives rise 
to two singular points of the angle $\pi$. We denote these points by  
$(0,0_+)$ and $(0,0_-)$. The action of $\gamma_{v}\gamma_{h}^{-1}\gamma_v$ 
on $T$ can be lifted to  the action  on the folded plane in
such a way that 
point $(0,0_+)$ is fixed. All singular points of type $(m,3n_+)$ can be 
naturally identified with $H_1(T,\mathbb{Z})$. So these points are 
invariant under the fourth power of the automorphism 
$\gamma_{v}\gamma_{h}^{-1}\gamma_v$ (lemma 3). 

Denote by $\tilde{\mathcal F_1}$ and $\tilde{\mathcal F_2}$ the
liftings of the foliations $\mathcal F_1$ and $\mathcal F_2$ 
on  the universal covering of the torus $T$.
Denote by $\tilde{f_1}$ the leaf of $\tilde{\mathcal F_1}$ starting at
the point $(0,0_+)$ and by $\tilde{f_2}$ the leaf of $\tilde{\mathcal F_2}$ 
starting in $(3,0_+)$

\begin{lemma}  The leaf $\tilde{f_1}$ of the foliation 
$\tilde{\mathcal F_1}$ comes arbitrary close to the point $(3,0_+)$. 
\end{lemma}

Proof. The point $(\frac{3}{2},\frac{\sqrt{3}}{2})$ is an intersection of
the leaves $\tilde{f_1}$ and $\tilde{f_2}$. The transformation 
$(\gamma_{v}\gamma_{h}^{-1}\gamma_v)^4$ 
preserves points $(0,0_+)$ and $(3,0_+)$ so it preserves leaves $\tilde{f_1}$ 
and $\tilde{f_2}$. It means that the point 
$(\gamma_{v}\gamma_{h}^{-1}\gamma_v)^4(\frac{3}{2},\frac{\sqrt{3}}{2})$ 
is also a point of the intersection 
of $\tilde{f_1}$ and $\tilde{f_2}$. Now use the fact that the action of 
$(\gamma_{v}\gamma_{h}^{-1}\gamma_v)^4$ contracts $\tilde{f_2}$. 
Applying the map $(\gamma_{v}\gamma_{h}^{-1}\gamma_v)^4$ inductively we 
construct a sequence of intersections of $\tilde{f_1}$ and
$\tilde{f_2}$ which tends to the 
point $(3,0_+)$. \hfill $\square$

Now we give a proof of the theorem 1 for the foliation 
$\tilde{\mathcal F_1}$ and its leaf
$\tilde{f_1}$. Since  $\gamma_{v}\gamma_{h}^{-1}\gamma_v$ performs a cyclic 
permutation of points $(\pm3,0_+)$, $(0,\pm 1_+)$ (lemma 1) the 
leaf $\tilde f_1$ approaches arbitrary close to all these points. 
It means that the closure of $\tilde f_1$ in $\mathbb{R}^2$ contains
the leaves of $\tilde{\mathcal F_1}$ passing through fore points
$(\pm3,0_+)$, $(0,\pm 1_+)$. Since $\tilde{\mathcal F_1}$ is periodic
in $\mathbb{R}^2$ the same argument tells that the leaf $\tilde{f_1}$ 
approaches arbitrary close to every singular point of the type $(3n,k_+)$. 
So its closure in the plane contains the closure of the union of 
all leaves of the foliation  $\tilde{\mathcal F_1}$ passing through 
singular points of type $(3n,k_+)$. The last closure is the 
whole plane (lemma 2). \hfill $\square$

{\bf Final remarks}. The leaf $\tilde{f_1}$ constructed in the theorem 1 is 
everywhere dense in the plane. One can see that this leaf 
propagates in the plane  with a logarithmic speed. 

Here is an optical interpretation of the constructed foliation. Consider the 
two dimensional plane $\mathbb{R}^2$ with the integer 
lattice $\mathbb{Z}^2$. Put in any vertex of the lattice a $90^o$ 
rotation invariant cross in such a way that all the picture is 
$\mathbb{Z}^2$ invariant (fig.3). Send a light in the direction parallel 
to the bissectrix of crosses and consider it reflections from 
edges of crosses. The constructed ray  corresponds to a leaf of some
foliation on the two-dimensional torus with 3 singular points. 

\begin{center}
\
\epsfbox{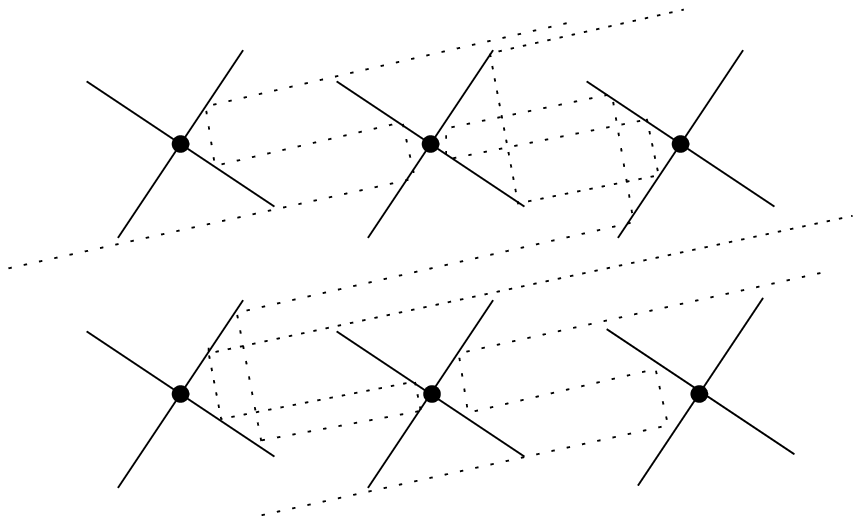}
\end{center}

{\bf Acknowledgments}. The author is grateful to M.Kontsevich for
introducing to the subject and very helpful discussions, to A.Glutsuk  
for communication of results of the article [An1] and to A.Zorich
for his interest to the work.

\end{document}